\documentclass[12pt]{amsart}
\usepackage{amssymb,amsmath}
\def\no#1{||#1||}
\def\nno#1{|||#1|||}
\numberwithin{equation}{section}
\def\db{\bar\partial}
\def\db*{\bar\partial^*}
\def\T{\text}

\def\simleq{\underset\sim<}
\def\1#1{\overline{#1}}
\def\2#1{\widetilde{#1}}
\def\3#1{\widehat{#1}}
\def\4#1{\mathbb{#1}}

\def\5#1{\frak{#1}}
\def\6#1{{\mathcal{#1}}}
\def\C{{\4C}}
\def\R{{\4R}}

\def\Z{{\4Z}}

\begin{document}
\abstract
We prove local hypoellipticity of the complex Laplacian $\Box$ in a domain  which has compactness estimates, is of finite type outside a curve transversal to the CR directions and for which the holomorphic tangential derivatives of a defining  function are subelliptic multipliers in the sense of Kohn. 
\newline
MSC: 32F10, 32F20, 32N15, 32T25 
\endabstract
\title[Hypoellipticity of the $\bar\partial$-Neumann problem...]{Hypoellipticity of the $\bar\partial$-Neumann problem at exponentially degenerate points}
\author[T.V.~Khanh and G.~Zampieri]{Tran Vu Khanh and  Giuseppe Zampieri}
\address{Dipartimento di Matematica, Universit\`a di Padova, via 
Trieste 63, 35121 Padova, Italy}
\email{khanh@math.unipd.it,  
zampieri@math.unipd.it}
\maketitle
\def\Giialpha{\mathcal G^{i,i\alpha}}
\def\cn{{\C^n}}
\def\cnn{{\C^{n'}}}
\def\ocn{\2{\C^n}}
\def\ocnn{\2{\C^{n'}}}
\def\const{{\rm const}}
\def\rk{{\rm rank\,}}
\def\id{{\sf id}}
\def\aut{{\sf aut}}
\def\Aut{{\sf Aut}}
\def\CR{{\rm CR}}
\def\GL{{\sf GL}}
\def\Re{{\sf Re}\,}
\def\Im{{\sf Im}\,}
\def\codim{{\rm codim}}
\def\crd{\dim_{{\rm CR}}}
\def\crc{{\rm codim_{CR}}}
\def\phi{\varphi}
\def\eps{\varepsilon}
\def\d{\partial}
\def\a{\alpha}
\def\b{\beta}
\def\g{\gamma}
\def\G{\Gamma}
\def\D{\Delta}
\def\Om{\Omega}
\def\k{\kappa}
\def\l{\lambda}
\def\L{\Lambda}
\def\z{{\bar z}}
\def\w{{\bar w}}
\def\Z{{\1Z}}
\def\t{{\tau}}
\def\th{\theta}
\emergencystretch15pt
\frenchspacing
\newtheorem{Thm}{Theorem}[section]
\newtheorem{Cor}[Thm]{Corollary}
\newtheorem{Pro}[Thm]{Proposition}
\newtheorem{Lem}[Thm]{Lemma}
\theoremstyle{definition}\newtheorem{Def}[Thm]{Definition}
\theoremstyle{remark}
\newtheorem{Rem}[Thm]{Remark}
\newtheorem{Exa}[Thm]{Example}
\newtheorem{Exs}[Thm]{Examples}
\def\Label#1{\label{#1}}
\def\bl{\begin{Lem}}
\def\el{\end{Lem}}
\def\bp{\begin{Pro}}
\def\ep{\end{Pro}}
\def\bt{\begin{Thm}}
\def\et{\end{Thm}}
\def\bc{\begin{Cor}}
\def\ec{\end{Cor}}
\def\bd{\begin{Def}}
\def\ed{\end{Def}}
\def\br{\begin{Rem}}
\def\er{\end{Rem}}
\def\be{\begin{Exa}}
\def\ee{\end{Exa}}
\def\bpf{\begin{proof}}
\def\epf{\end{proof}}
\def\ben{\begin{enumerate}}
\def\een{\end{enumerate}}
\def\dotgamma{\Gamma}
\def\dothatgamma{ {\hat\Gamma}}

\def\simto{\overset\sim\to\to}
\def\1alpha{[\frac1\alpha]}
\def\T{\text}
\def\R{{\Bbb R}}
\def\I{{\Bbb I}}
\def\C{{\Bbb C}}
\def\Z{{\Bbb Z}}
\def\Fialpha{{\mathcal F^{i,\alpha}}}
\def\Fiialpha{{\mathcal F^{i,i\alpha}}}
\def\Figamma{{\mathcal F^{i,\gamma}}}
\def\Real{\Re}
%
%
%
\section{  Introduction}
For the pseudoconvex domain $\Omega\subset\C^n$ whose boundary is   defined in coordinates $z=x+iy$ of $\C^n$,  by
\begin{equation}
\Label{1.1}
2x_n=\exp\left({-\frac1{(\sum_{j=1}^{n-1}|z_j|^2)^{\frac s2}}}\right),\quad s>0,
\end{equation}
the tangential Kohn Laplacian $\Box_b =\bar\partial_b\bar\partial^{ *}_b+\bar\partial^{*}_b\bar\partial_b$ as well as the full Laplacian $\Box=\bar\partial\bar\partial^{ *}+\bar\partial^{*}\bar\partial$ show very interesting features  especially  in comparison with the ``tube domain" whose boundary is defined by
\begin{equation}
\Label{1.2}
2x_n=\exp\left({-\frac1{(\sum_{j=1}^{n-1}|x_j|^2)^{\frac s2}}}\right),\quad s>0.
\end{equation}
(Here $z_j$ have been  replaced by $x_j$ at exponent.) Energy estimates are the same for the two domains. For the  problem on the boundary $b\Omega$, they come as
\begin{multline}
\Label{1.2bis}
\no{(\log \Lambda )^{\frac1s}u}_{b\Omega}\simleq \no{\bar\partial_b u}^2_{b\Omega}+\no{\bar\partial^{*}_b u}^2_{b\Omega}+\no{u}^2_{b\Omega}\\
\T{for any smooth compact support  form $u\in C^\infty_c(b\Omega)^k$ of degree $k\in [1,n-2]$}.
\end{multline}
Here $\log \Lambda$ is the tangential pseudodifferential operator with symbol $\log(1+|\xi'|^2)^{\frac12}),\,\xi'\in\R^{2n-1}$, the dual real tangent space. As for the  problem on the domain $\Omega$, one has simply to replace $\bar\partial_b,\,\bar\partial^*_b$ by $\bar\partial,\,\bar\partial^*$  and take norms over $\Omega$ for forms $u$ in $D_{\bar\partial^*}$, the domain of $\bar\partial^*$, of degree $1\le k\le n-1$; this can be seen, for instance, in   \cite{KZflat09}.
In particular, these are superlogarithmic (resp. compactness) estimates if $s<1$ (resp. for any $s>0$). A related problem is that of the local hypoellipticity of  the Kohn Laplacian $\Box_b $ or, with equivalent terminology, the local  regularity  of the inverse (modulo harmonics)  operator $N_b=\Box^{-1}_b$. Similar is the notion of  hypoellipticity of the Laplacian $\Box$ or the regularity of the inverse Neumann operator $N=\Box^{-1}$. 
It has been proved by Kohn in \cite{K02} that superlogarithmic estimates suffice for local hypoellipticity  of the problem both in the boundary and in the domain.
(Note that hypoellipticity for the domain, \cite{K02} Theorem~8.3, is deduced from microlocal hypoellipticity for the boundary, \cite{K02} Theorem~7.1, but a direct proof is also available, \cite{Kh09} Theorem~5.4.)
 In particular, for  \eqref{1.1} and \eqref{1.2}, there is local hypoellipticity when $s<1$. 

As for the more delicate  hypoellipticity, in the uncertain range of indices $s\ge 1$, only the tangential problem  has been studied and the striking conclusion is that  the behavior of \eqref{1.1} and \eqref{1.2} split. The first stays always hypoelliptic for any $s$ (Kohn \cite{K00}) whereas the second is not for $s\ge 1$ (Christ \cite{Ch02}). When one tries to relate $(\bar\partial_b,\bar\partial^{*}_b)$ on $b\Omega$ to $(\bar\partial,\bar\partial^*) $ on $\Omega$, estimates go well through (Kohn \cite{K02} Section 8 and Khanh \cite{Kh09} Chapter 4) but not regularity. In particular, the two conclusions about tangential hypoellipticity of $\Box_b $ for \eqref{1.1} and non-hypoellipticity for \eqref{1.2} when $s\ge1$, cannot be automatically transferred from $b\Omega$ to $\Omega$. Now, for the non-hypoellipticity in $\Omega$ in case of the tube \eqref{1.2} we have obtained with Baracco in \cite{BKZ09} a result of propagation which is not  equivalent but intimately related. The real lines $x_j$ are propagators of holomorphic extendibility from $\Omega$ across $b\Omega$. 
What we prove in the present paper is  hypoellipticity in $\Omega$ for \eqref{1.1} when $s\geq 1$.
\bt
\Label{t1.1}
Let $\Omega$ be a pseudoconvex domain of $\C^n$   in a neighborhood of $z_o=0$  and assume that the $\bar\partial$-Neumann problem satisfies the following properties
\begin{itemize}
\item[(i)] there are local compactness estimates,
\item[(ii)] there are subelliptic estimates for $(z_1,...,z_{n-1})\neq0$,
\item[(iii)] $\partial_{z_j}r,\,\,j=1,...,n-1$, are subelliptic multipliers (cf. \cite{K79}).
\end{itemize}
Then $\Box$ is  locally hypoelliptic at $z_o$.
\et
The proof follows in Section 2. It consists in relating the system on $\Omega$ to the tangential system on $b\Omega$  along the guidelines of \cite{K02} Section 8, and then in using the argument of \cite{K00} simplified by the additional assumption $(i)$. 
\br
The domain with boundary \eqref{1.1}, but not \eqref{1.2}, satisfies the hypotheses of Theorem~\ref{t1.1} for any $s>0$: (i) is obvious, and (ii) and (iii) are the content of \cite{K00} Section 4.
\er
Notice that $\partial\Omega $ is given only locally in a neighborhodd of $z_o$. We can continue $\partial\Omega$ leaving it unchanged in a neighborhood of $z_o$, making it strongly pseudoconvex elsewhere, in such a way that it bounds a relatively compact domain $\Omega\subset\subset\C^n$ (cf. \cite{MN92}).
In this situation $\Box$ is hypoelliptic at every boundary point. Also, it is well defined a $H^0$ inverse Neumann operator $N=\Box^{-1}$, and, by Theorem~\ref{t1.1}, the $\bar\partial$-Neumann solution operator $\bar\partial^* N$ preserves $C^\infty(\bar\Omega)$-smoothness. It even preserves the exact Sobolev class $H^s$ according to Theorem~\ref{t2.2} below.  In other words, the canonical solution $u=\bar\partial^*Nf$ of $\bar\partial u=f$ for $f\in\T{Ker }\bar\partial$ is $H^s$ exactly at the points of $b\Omega$ where $f$ is $H^s$. 
The Bergman 
projection $B$ also preserves $C^\infty(\bar\Omega)$-smoothness 
 on account of Kohn's formula $B=\T{Id}-\bar\partial^*N\bar\partial$.
 \vskip0.5cm
 \noindent
{\it Aknowledgments. } The authors are grateful to Emil Straube for suggesting the argument which leads to the hypoellipticity of the operator $\Box$ from that of the system $(\bar\partial,\bar\partial^*,\Delta)$.
\section{Hypoellipticity of $\Box$ and exact hypoellipticity of $\bar\partial^*N$}
We state properly  hypoellipticity and exact hypoellipticity of a general system $(P_j)$.
\bd
\Label{d2.1}
(i) The system $(P_j)$  is locally hypoelliptic  at $z_o\in b\Omega$ if 
$$
P_ju\in C^\infty(\bar\Omega)_{z_o}^k\T{ for any $j$ implies $u\in C^\infty(\bar\Omega)_{z_o}^k$},
$$
where $C^\infty(\bar\Omega)_{z_o}^k$ denotes the set of germs of $k$-forms smooth at $z_o$.

\noindent
(ii)
The system $(P_j)$  is exactly locally hypoelliptic  at  $z_o\in b\Omega$ when there is a neighborhood $U$ of $z_o$ such that for any pair of cut-off functions $\zeta$ and $\zeta'$ in $C^\infty_c(U)$ with $\zeta'|_{\T{supp}(\zeta)}\equiv1$  we have for any $s$ and for suitable $c_s$
\begin{equation}
\Label{2.-1}
\no{\zeta u}^2_s\leq c_s(\sum_j\no{\zeta'P_ju}^2_s+\no{u}^2_0),\qquad u\in C^\infty(\bar\Omega)^k\cap D_{(P_j)}.
\end{equation}
\ed
If $(P_j)$ happens to have an inverse, this is said to be locally regular and locally exactly regular in the situation of (i) and (ii) respectively.
\br
By Kohn-Nirenberg \cite{KN65} the  assumption $u\in C^\infty$ can be removed from \eqref{2.-1}. Precisely, by the elliptic regularization, one can prove that if $\zeta'P_j u\in H^s$ and $\zeta' u\in H^0$, then $\zeta u\in H^s$ and satisfies \eqref{2.-1}. This motivates the word ``exact", that is, Sobolev exact. Not only the local $C^\infty$- but also the $H^s$-smoothness passes from $P_ju$ to $u$. 
\er

Let $\vartheta$ be the formal adjoint of $\bar\partial$ and $\Delta=\bar\partial\vartheta+\vartheta\bar\partial$  the Laplacian; it acts on forms by the action of the usual Laplacian on its coefficients. If $u\in D_\Box$, then $\Box u=\Delta u$. We first prove exact hypoellipticity of the system $(\bar\partial,\bar\partial^*,\Delta)$; hypoellipticity of $\Box$ itself will follow by the method of Boas-Straube.
\bt
\Label{t2.1}
In the situation of Theorem~\ref{t1.1}, we have, for a neighborhood $U$ of $z_o$ and for any couple of cut-off $\zeta$ and $\zeta'$ with $\zeta'|{\T{supp}\,\zeta}\equiv1$
\begin{equation}
\Label{2.1}
\no{\zeta u}^2_s\simleq \no{\zeta'\bar\partial u}^2_s +\no{\zeta'\bar\partial^*u}^2_s+\no{\zeta'\Delta u}^2_{s-2}+\no{u}^2_0,\quad u\in D_{\bar\partial^*}.
\end{equation}
In particular, the system $(\bar\partial,\bar\partial^*,\Delta)$ is exactly locally hypoelliptic at $z_o=0$. 
\et
\br
The hypoellipticity of $\Box_b $ under (ii) and (iii) of Theorem~\ref{t1.1} is proved by Kohn in \cite{K00}. It does not require (i) but it is not exact hypoellipticity (the neighborhood $U$ of \eqref{2.-1} depends on $s$). However, inspection of his proof shows that, if (i) is added, then in fact \eqref{2.-1} holds for $(P_j)=\Box_b $. Our proof consists in a reduction to the tangential system.
\er
\bpf
We proceed in several steps which are highlighted in two intermediate propositions. 
We use the standard notation $Q(u,u)$ for $\no{\bar\partial u}^2_0+\no{\bar\partial^*u}^2_0$ and some variants as, for an operator $Op$,  $Q_{Op}(u,u):=\no{Op\, \bar\partial u}^2_0+\no{Op\,\bar\partial^*u}^2_0$; most often, in our paper, Op is chosen as $\Lambda^s\zeta'$.
We decompose a form $u$ as
\begin{equation*}
\begin{cases}
u=u^\tau+u^\nu,
\\
u^{\tau}=u^{\tau\,+}+u^{\tau\,-}+u^{\tau\,0},
\end{cases}
\end{equation*}
where the first is the decomposition in tangential and normal component and the second is the microlocal decomposition $u^{\tau\,\overset\pm0}=
\Psi^{\overset\pm0}u^\tau$ in which $\Psi^{\overset\pm0}$ are the tangential pseudodifferential operators whose symbols $\psi^{\overset\pm0}$ 
are a conic decomposition of the unity in the space dual to $\R^{2n-1}$ the real orthogonal to $\partial r$ (cf. Kohn \cite{K02}).  We begin our proof by remarking that
 any of the forms  $u^\#=u^\nu,
\,u^{\tau\,-},\,u^{\tau\,0}$ enjoys elliptic estimates
\begin{equation}
\Label{2.3}
\no{\zeta u^\#}_s^2\simleq\no{\zeta'\bar\partial u^\#}_{s-1}^2+\no{\zeta'\bar\partial^* u^\#}_{s-1}^2+\no{u^\#}^2_0\qquad s\geq2.
\end{equation}
We refer to \cite{FK72} formula (1) of Main theorem as a general reference but also give an outline of the proof. 
For this, we have to call into play the tangential $s$-Sobolev norm which is defined by $\nno{u}_s=\no{\Lambda^su}_0$. 
We start from
\begin{equation}
\Label{2.3bis}
\nno{\zeta u^\#}^2_1\simleq Q(\zeta u^\#,\zeta u^\#)+\no{u^\#}^2_0;
\end{equation}
this is the basic estimate for $u^\nu$ (which vanishes at $b\Omega$) whereas it is \cite{K02} Lemma 8.6 for $u^{\tau\,-}$ and $u^{\tau\,\,0}$. Applying 
\eqref{2.3bis} to $\zeta'\Lambda^{s-1}\zeta u^\#$ one gets the estimate of tangential norms for any $s$. Finally, by non-characteristicity of $(\bar\partial,\bar\partial^*)$ one passes from tangential to full norms along the guidelines of \cite{Z08} Theorem~1.9.7. The version of this argument for $\Box$ can be found in \cite{K02} second part of p. 245.   Because of \eqref{2.3}, it suffices to prove \eqref{2.1} for the only $u^{\tau\,+}$. 
We further decompose
$$
u^{\tau\,+}=u^{\tau\,+\,(h)}+u^{\tau\,+\,(0)},
$$
where $u^{\tau\,+\,(h)}$ is the ``harmonic extension" in the sense of Kohn \cite{K02} and $u^{\tau\,+\,(0)}$ is just the complementary part.
We denote by $\bar\partial^\tau$ the extension of $\bar\partial_b$ from $b\Omega$ to $\Omega$ which stays tangential to the level surfaces $r\equiv\T{const}$. It acts on tangential forms $u^\tau$ and it is defined by $\bar\partial^\tau u^\tau=(\bar\partial u^\tau)^\tau$. We denote by $\bar\partial^{\tau\,*}$ its adjoint; thus $\bar\partial^{\tau\,*}u^\tau=\bar\partial^*(u^\tau)$. We use the notations $\Box^\tau$ and $Q^\tau$ for the corresponding Laplacian and energy.
 We notice that over a tangential form $u^\tau$ we have a decomposition
\begin{equation} 
\Label{2.3ter}
Q={Q^\tau} +\no{\bar L_nu^{\tau}}^2_0.
\end{equation}
The proof of \eqref{2.1} for $u^{\tau\,+}$ requires two crucial technical results. Here is the first which is the most central
\bp
\Label{p2.1}
For the harmonic extension $u^{\tau\,+\,(h)}$ we have 
\begin{equation}
\Label{2.6}
\nno{\zeta u^{\tau\,+\,(h)}}^2_s\simleq Q^\tau_{\Lambda^s\zeta'}(u^{\tau\,+\,(h)},u^{\tau\,+\,(h)})+\no{u^{\tau\,+\,(h)}}^2_0.
\end{equation}
\ep
\bpf
We apply compactness estimates (cf. e.g. \cite{Kh09} Section 6) for $\zeta'\Lambda^s\zeta u^{\tau\,+\,(h)}$,
\begin{equation}
\Label{2.7}
\no{\zeta'\Lambda^s\zeta u^{\tau\,+\,(h)}}^2\leq\epsilon Q(\zeta'\Lambda^s\zeta u^{\tau\,+\,(h)},\zeta'\Lambda^s\zeta u^{\tau\,+\,(h)})+c_\epsilon\no{\zeta'\Lambda^s\zeta u^{\tau\,+\,(h)}}^2_{-1}.
\end{equation}
We decompose $Q$ according to \eqref{2.3ter}.
We calculate ${Q^\tau} $ over $\zeta'\Lambda^s\zeta u^{\tau\,+\,(h)}$ and compute errors coming from commutators $[{Q^\tau} ,\zeta'\Lambda^s\zeta]$. In this calculation we assume that the cut off functions are of product type  $\zeta(z')\zeta(t)$ where $z'$ (resp. $t$) are complex (resp. totally real) tangential coordinates in $T_{z_o}b\Omega$. We have
\begin{equation}
\Label{2.9}
\begin{split}
{Q^\tau}(\zeta'&\Lambda^s\zeta u^{\tau\,+\,(h)},\zeta'\Lambda^s\zeta u^{\tau\,+\,(h)})
\\
&\simleq {Q^\tau}_{\zeta'\Lambda^s\zeta}(u^{\tau\,+\,(h)},u^{\tau\,+\,(h)})+\nno{\zeta u^{\tau\,+\,(h)}}^2_s+\nno{\zeta'u^{\tau\,+\,(h)}}^2_{s-1}
\\
&+\left(\no{(|\dot \zeta(z')|+|\dot \zeta'(z')|)\Lambda^s u^{\tau\,+\,(h)}}^2_0+\no{\sum_{j=1}^{n-1}|r_{z_j}|(|\dot\zeta(t)|+|\dot\zeta'(t)|)\Lambda^s u^{\tau\,+\,(h)}}^2_0\right).
\end{split}
\end{equation}
We explain \eqref{2.9}. First, the commutators $[{\bar\partial^\tau},\zeta'\Lambda^s\zeta]$ (and similarly as for $[\bar\partial^{\prime *},\zeta'\Lambda^s\zeta]$) are decomposed by Jacobi identity as
$$
[{\bar\partial^\tau},\zeta'\Lambda^s\zeta]=[{\bar\partial^\tau},\zeta']\Lambda^s\zeta+\zeta'[{\bar\partial^\tau},\Lambda^s]\zeta+\zeta'\Lambda^s[{\bar\partial^\tau},\zeta].
$$
The central commutator $[{\bar\partial^\tau},\Lambda^s]$ produces the error term $\nno{\zeta u^{\tau\,+\,(h)}}^2_s$. As for the two others, we have
$$
[{\bar\partial^\tau},\zeta(z')\zeta(t)]=[{\bar\partial^\tau},\zeta(z')]\zeta(t)+\zeta(z')[{\bar\partial^\tau},\zeta(t)],
$$ 
and similarly for $\zeta $ replaced by $\zeta'$ and ${\bar\partial^\tau}$ by ${\bar\partial^{\tau\,*}}$. Now,
\begin{equation}
\Label{2.10}
[{\bar\partial^\tau},\zeta(z')]\sim\dot\zeta(z').
\end{equation}
On the other hand, 
we first notice that it is not restrictive to assume that $\partial_{z_1},...,\partial_{z_{n-1}}$ are a basis of $T^{1,0}_0b\Omega$ for otherwise, owing to (iii), we have subelliptic estimates from which local regularity readily follows. Thus, 
 each $\bar L_j,\,\,j=1,...,n-1$, is of type $\bar L_j=r_{\bar z_j}\partial_{\bar z_n}-r_{\bar z_n}\partial_{\bar z_j}$, and then
\begin{equation}
\Label{2.11}
\begin{split}
[{\bar\partial^\tau},\zeta(t)]&\sim\sum_{j=1}^{n-1}[\bar L_j,\zeta(t)]
\\
&\sim\sum_{j=1}^{n-1} r_{\bar z_j}\dot\zeta(t).
\end{split}
\end{equation}
By combining \eqref{2.10} with \eqref{2.11} (and using the analogous for $\zeta'$ and ${\bar\partial^{\tau\,*}}$), we get the last line of \eqref{2.9}. This establishes \eqref{2.9}. Next, since $({\bar\partial^\tau},{\bar\partial^{\tau\,*}})$ has subelliptic estimates, say $\eta$-subelliptic, for $z'\neq0$ and hence in particular over $\T{supp}\,\dot\zeta(z')$ and $\T{supp}\,\dot\zeta'(z')$ and since the $r_{\bar z_j}$ are, say, $\eta$-subelliptic multipliers even at $z'=0$, then the last line of \eqref{2.9} is estimated by $\no{\zeta''\Lambda^{s-\eta}\zeta'u^{\tau\,+\,(h)}}^2$ where $\zeta''\equiv1$ over $\T{supp}\,\zeta'$. This shows, using iteration over increasing $k$ such that $k\eta>s$ and over decreasing $j$ from $s-1$ to $0$, that \eqref{2.7} and \eqref{2.9} imply \eqref{2.6} provided that we add on the right side the extra  term $\no{\bar L_n\zeta'\Lambda^s\zeta u^{\tau\,+\,(h)}}^2$. Note that, as a result of the inductive process, we have to replace $Q_{\zeta'\Lambda^s\zeta}$ in \eqref{2.9} by $Q_{\Lambda^s\zeta'}$ in \eqref{2.6}.

Up to this point the argument is the same as in \cite{K00} and does not make any use of the specific properties of the harmonic extension  $u^{\tau\,+\,(h)}$. We start the new part which  is dedicated to prove that $\no{\bar L_n\zeta'\Lambda^s\zeta u^{\tau\,+\,(h)}}^2$ can be removed from the right of \eqref{2.6}. For this we have to use the main property of this extension expressed by \cite{K02} Lemma 8.5, that is,
\begin{equation}
\Label{2.12}
\no{\bar L_n \zeta u^{\tau\,+\,(h)}}^2_0\simleq \sum_{j=1}^{n-1}\no{\bar L_j \zeta u_b^{\tau\,+}}^2_{b,\,-\frac12}+\no{u^{\tau\,+}}^2_0.
\end{equation}
Note that \eqref{2.12} differs from \cite{K02} Lemma 8.5 by $[\bar L_n,\Psi^+]$; but this is an error term 
which can be taken care of by $u^{\tau\,0}$ to which elliptic estimates apply.
Applying \eqref{2.12} to $\zeta'\Lambda^s\zeta u^{\tau\,+\,(h)}$ (for the first inequality below), and  using the classical inequality $\no{\cdot}^2_{b,\,-\frac12}\leq c_\epsilon\no{\cdot}^2_0+\epsilon\nno{\partial_r\cdot}^2_{-1}$ (cf. e.g. \cite{KZtan09} (1.10)) together with the splitting $\partial_r=\bar L_n+Tan$ (for the second), we get 
\begin{equation}
\Label{2.13}
\begin{split}
\no{\bar L_n\zeta'\Lambda^s\zeta u^{\tau\,+\,(h)}}^2_0&\underset{\T{by \eqref{2.12}}}\simleq \sum_{j=1}^{n-1}\no{\bar L_j\zeta'\Lambda^s\zeta u^{\tau\,+}_b}^2_{b,\,-\frac12}+\no{\zeta'\Lambda^s\zeta u^{\tau\,+}}^2_0
\\
&\simleq c_{\epsilon}\sum_{j=1}^{n-1}\no{\bar L_j\zeta'\Lambda^s\zeta u^{\tau\,+\,(h)}}^2_0+\epsilon\sum_{j=1}^{n-1}\nno{\bar L_n\bar L_j\zeta'\Lambda^s\zeta u^{\tau\,+\,(h)}}^2_{-1}
\\
&+\epsilon\sum_{j=1}^{n-1}\nno{Tan \,\bar L_j\zeta'\Lambda^s\zeta u^{\tau\,+\,(h)}}_{-1}^2+\no{\zeta'\Lambda^s\zeta u^{\tau\,+\,(h)}}^2_0.
\end{split}
\end{equation}
The first term on the right of the last inequality is controlled by $\underset{j=1}{\overset{n-1}\sum}\no{\zeta'\Lambda^s\zeta \bar L_j u^{\tau\,+\,(h)}}^2+\nno{\zeta u^{\tau\,+\,(h)}}^2_s+\nno{\zeta'' u^{\tau\,+\,(h)}}^2_{s-1}$ by the first part of the proposition; moreover, we have the immediate estimate $\sum_{j=1}^{n-1}\no{\zeta'\Lambda^s\zeta \bar L_j u^{\tau\,+\,(h)}}^2\simleq Q^\tau_{\Lambda^s\zeta'}(u^{\tau\,+\,(h)},u^{\tau\,+\,(h)})$.
The term which carries $\epsilon\,Tan$, after $Tan$ has been annihilated by the Sobolev norm of index $-1$, has the same estimate as the first term. It remains to control the second term in the right which involves $\epsilon\bar L_n$. We rewrite $\bar L_n\bar L_j=\bar L_j\bar L_n+[\bar L_n,\bar L_j]$; when $\bar L_j$ moves in first position, it is annihilated by $-1$ and what remains is absorbed in the left. As for the commutator, we have
\begin{equation*}
\begin{split}
\nno{[\bar L_n,\bar L_j]\zeta'\Lambda^s\zeta u^{\tau\,+\,(h)}}^2_{-1}&\simleq \nno{\zeta u^{\tau\,+\,(h)}}^2_s+\nno{\partial_r\zeta'\Lambda^s\zeta u^{\tau\,+\,(h)}}^2_{-1}
\\
&\simleq \nno{\zeta u^{\tau\,+\,(h)}}^2_s+\nno{\bar L_n\zeta'\Lambda^s\zeta u^{\tau\,+\,(h)}}^2_{-1},
\end{split}
\end{equation*}
where we have used the splitting $\partial_r=Tan +\bar L_n$ in the second inequality. Again, the term with $\bar L_n$, which now comes in $-1$ norm, is absorbed in the left of \eqref{2.13}.
Summarizing up, we have got
\begin{equation}
\Label{2.14}
\begin{split}
\no{\bar L_n\zeta'\Lambda^s\zeta u^{\tau\,+\,(h)}}^2_0&\simleq c_\epsilon Q^\tau_{\Lambda^s\zeta'}(u^{\tau\,+\,(h)},u^{\tau\,+\,(h)})
\\
&+\nno{\zeta u^{\tau\,+\,(h)}}^2_s+\nno{\zeta'' u^{\tau\,+\,(h)}}^2_{s-1}.
\end{split}
\end{equation}
But $\nno{\bar L_n\cdot}^2$ comes with a factor $\epsilon$ of compactness and hence the term in $s$-norm in the last line can be absorbed in the left of the initial inequalities \eqref{2.7} or \eqref{2.6}.
Finally, we use an inductive argument to go down from $s-1$ to $0$.
 This concludes the proof of the proposition.

\epf
We remark now that
\begin{equation}
\Label{2.18}
\begin{split}
\no{\zeta u^{\tau\,+\,(h)}}_0^2&\simleq \no{\zeta u_b^{\tau\,+}}^2_{b,\,-\frac12}
\\
&\simleq \no{\zeta u^{\tau\,+}}^2_0+\nno{\partial_r\zeta u^{\tau\,+}}_{-1}^2
\\
&\leq \no{\zeta u^{\tau\,+}}^2_0+\nno{\bar L_n\zeta u^{\tau\,+}}^2_{-1}+\nno{Tan\, \zeta u^{\tau\,+}}^2_{-1}
\\
&\simleq Q_{\Lambda^{-1}\zeta}(u^{\tau\,+},u^{\tau\,+})+\no{\zeta u^{\tau\,+}}^2_0.
\end{split}
\end{equation}
The same inequality also holds for $u^{\tau\,+\,(h)}$ replaced by $u^{\tau\,+\,(0)}$ on account of the identity $u^{\tau\,+\,(0)}=u^{\tau\,+}+u^{\tau\,+\,(h)}$.
We need another preparation result
\bp
\Label{p2.2}
We have
\begin{equation}
\Label{2.16}
{Q^\tau}_{\Lambda^s\zeta'}(u^{\tau\,+\,(h)},u^{\tau\,+\,(h)})\simleq {Q^\tau}_{\Lambda^s\zeta'}(u^{\tau\,+},u^{\tau\,+})+{Q^\tau}_{\partial_r\Lambda^{s-1}\zeta'}(u^{\tau\,+},u^{\tau\,+})
\end{equation}
and
\begin{equation}
\Label{2.15}
\begin{split}
\nno{\zeta u^{\tau\,+\,(0)}}^2_s&\simleq {Q^\tau}_{\Lambda^{s-1}\zeta'}(u^{\tau\,+},u^{\tau\,+})+{Q^\tau}_{\partial_r\Lambda^{s-2}\zeta'}(u^{\tau\,+},u^{\tau\,+})
\\
&\quad+\nno{\zeta' \Delta u^{\tau\,+}}^2_{s-2}
+\no{ u^{\tau\,+}}_{0}^2.
\end{split}
\end{equation}
\ep
\bpf
The proof of \eqref{2.16} is an immediate combination of the formulas $\no{\zeta'u^{\tau\,+\,(h)}}_0\simleq \no{\zeta'u^{\tau\,+}_b}_{b,\,-\frac12}$ and $\no{\zeta'u^{\tau\,+}}_{b,\,-\frac12}\simleq \no{\zeta'u^{\tau\,+}}_0+\nno{\partial_r \zeta'u^{\tau\,+}}^2_{-1}$. 

We prove now \eqref{2.15}. By elliptic estimate for $u^{\tau\,+\,(0)}$ (which vanishes at $b\Omega$) with respect to the order 2 elliptic operator $\Delta$, we have
\begin{equation}
\Label{2.17}
\nno{\zeta u^{\tau\,+\,(0)}}_s^2\simleq \nno{\zeta'\Delta u^{\tau\,+\,(0)}}^2_{s-2}+\no{u^{\tau\,+\,(0)}}^2_{0}.
\end{equation}
This result of Sobolev regularity at the boundary is very classical: it is formulated, for functions in $H^1_0$ such as the coefficients of $u^{\tau\,+\,(0)}$, e.g. in Evans \cite{E97} Theorem 5 p. 323.
Owing to the identity  $\Delta u^{\tau\,+\,(0)}=\Delta u^{\tau\,+}+P^1u^{\tau\,+\,(h)}$ for a 1-order operator $P^1$ (cf. \cite{K02} p. 241), we can replace $\Delta u^{\tau\,+\,(0)}$ by $\Delta u^{\tau\,+}$ on the right side of \eqref{2.17} putting the contribution of $P^1$ into an error term of type $\nno{\zeta' u^{\tau\,+\,(h)}}_{s-1}+\nno{\zeta'\partial_ru^{\tau\,+\,(h)}}_{s-2}$,
which can be estimated, on account of  the splitting $\partial_r=\bar L_n+Tan$, by $\nno{\zeta'u^{\tau\,+\,(h)}}_{s-1}+\nno{\zeta''u^{\tau\,+\,(h)}}_{s-2}+{Q^\tau}_{\Lambda^{s-2}\zeta'}(u^{\tau\,+\,(h)},u^{\tau\,+\,(h)})$. 
We write the terms of order $s-1$ and $s-2$ as a common $\nno{\zeta''u^{\tau\,+\,(h)}}_{s-1}$ that we can estimate, using \eqref{2.6} and \eqref{2.16}, by
$$
\nno{\zeta'' u^{\tau\,+\,(h)}}_{s-1}^2\simleq Q^\tau_{\Lambda^{s-1}\zeta'''}(u^{\tau\,+},u^{\tau\,+})+Q^\tau_{\Lambda^{s-2}\partial_r\zeta'''}(u^{\tau\,+
},u^{\tau\,+}).
$$
This brings down from $s-1$ to $0$ the Sobolev index in the error term. This $0$-order term  $\no{u^{\tau\,+\,(h)}}^2_0$, together with its companion  $\no{u^{\tau\,+\,(0)}}^2_0$  in the right of \eqref{2.17}, is estimated, because of \eqref{2.18}, by $\no{u^{\tau\,+}}^2_0$ up to a term 
$Q_{\Lambda^{-1}\zeta}$ which is controlled by the right side of \eqref{2.15}. This concludes the proof of \eqref{2.15}.

\epf

\vskip0.4cm
\noindent
{\it End of proof of Theorem~\ref{t2.1}. }
We  prove \eqref{2.1} for $u^{\tau\,+}$; this implies the conclusion in full generality according to the first part of the proof. We have
\begin{equation}
\Label{2.22}
\begin{split}
\nno{\zeta&u^{\tau\,+\,(h)}}^2_s\underset{\T{by \eqref{2.6}}}\simleq {Q^\tau}_{\Lambda^s\zeta'}(u^{\tau\,+\,(h)},u^{\tau\,+\,(h)})+\no{u^{\tau\,+\,(h)}}^2_0
\\
&\underset{\T{by \eqref{2.16} and \eqref{2.18}}}\simleq {Q^\tau}_{\Lambda^s\zeta'}(u^{\tau\,+},u^{\tau\,+})+ {Q^\tau}_{\partial_r\Lambda^{s-1}\zeta'}(u^{\tau\,+},u^{\tau\,+})+\no{u^{\tau\,+}}^2_0.
\end{split}
\end{equation}
We combine \eqref{2.22} with \eqref{2.15}; what we get is
\begin{equation}
\Label{2.20}
\begin{split}
\nno{\zeta u^{\tau\,+}}^2_s&\leq \nno{\zeta u^{\tau\,+\,(h)}}^2_s+\nno{\zeta u^{\tau\,+\,(0)}}^2_s
\\
&\simleq \no{\zeta'\bar\partial u^{\tau\,+}}^2_s+\no{\zeta'\bar\partial^* u^{\tau\,+}}^2_s+\nno{\zeta'\Delta u^{\tau\,+}}^2_{s-2}+\no{u^{\tau\,+}}^2_0.
\end{split}
\end{equation}
By the non-characteristicity of $Q$, we can replace the tangential norm $\nno{\cdot}_s$ by the full norm $\no{\cdot}_s$ in the left of \eqref{2.20}.
(The explanation of this point can be found, for example, in \cite{K02} second part of p. 245.)
This proves \eqref{2.1} for $u^{\tau\,+}$ and thus also for a general $u$. 

\epf
We modify $b\Omega$ outside a neighborhood of $z_o$ where it satisfies the hypotheses of Theorem~\ref{t1.1} so that it is strongly pseudoconvex in the 
modified portion and bounds a relatively compact domain; in particular, there is well defined the $H^0$ inverse $N$ of $\Box$ in this domain.
There is an immediate crucial consequence of Theorem~\ref{t2.1}.
\bt
\Label{t2.2}
We have that 
\begin{equation}
\Label{2.30}
\bar\partial^*N\T{ is exactly regular over }\T{Ker}\,\bar\partial
\end{equation}
and
\begin{equation}
\Label{2.31}
\bar\partial N\T{ is exactly regular over }\T{Ker}\,\bar\partial^*.
\end{equation}
\et
\bpf
As for \eqref{2.30}, we put $u=\bar\partial^*Nf$ for $f\in\T{Ker}\,\bar\partial$. We get 
\begin{equation*}
\begin{cases}
\bar\partial u=f,
\\
\bar\partial^*u=0,
\\
\begin{split}
\Delta u&=(\vartheta\bar\partial+\bar\partial\vartheta)\bar\partial^*Nf
\\
&=\vartheta(\bar\partial\bar\partial^*+\bar\partial^*\bar\partial)Nf+\bar\partial\vartheta\db*Nf
\\&=\vartheta\Box Nf=\vartheta f.
\end{split}
\end{cases}
\end{equation*}
Thus, by \eqref{2.1}
\begin{equation}
\Label{2.32}
\begin{split}
\no{\zeta u}^2_s&\simleq \no{\zeta' f}^2_s+\no{\zeta'\vartheta f}^2_{s-2}+\no{u}^2_0
\\&\simleq \no{\zeta'f}^2_s+\no{u}^2_0.
\end{split}
\end{equation}
To prove \eqref{2.31}, we put $u=\bar\partial  Nf$ for $f\in\T{Ker}\,\db*$. We have a similar calculation as above which leads to the same formula as \eqref{2.32} (with the only difference that $\vartheta$ is replaced by $\bar\partial $ in the intermediate inequality). Thus from \eqref{2.32} applied both for $\bar\partial^* N$ and $\bar\partial N$ on $\T{Ker}\, \bar\partial $ and $\T{Ker}\, \db*$ respectively, we conclude  that these operators are exactly regular.

\epf
We are ready for the proof of Theorem~\ref{t1.1}. This follows from Theorem~\ref{t2.2} by the method of Boas-Straube.
\vskip0.3cm
\noindent
{\it Proof of Theorem~\ref{t1.1}.}
From the regularity of $\bar\partial^* N$ it follows that the Bergman projection $B$ is also regular. (Notice that exact regularity is perhaps lost by taking $\bar\partial$ in $B$.) 
We exploit formula (5.36) in \cite{S10} in unweighted norms, that is, for $t=0$:
\begin{equation*}
\begin{split}
N_q&=B_q(N_q\bar\partial)(\T{Id}-B_{q-1})(\db*N_q)B_q
\\&\quad +(\T{Id}-B_q)(\db*N_{q+1})B_{q+1}(N_{q+1}\bar\partial )(\T{Id}-B_q).
\end{split}
\end{equation*}
Now, in the right side, the $\bar\partial  N$'s and $\db* N$'s are evaluated over $\T{Ker}\, \db*$ and $\T{Ker}\,\bar\partial $ respectively; thus they are exactly regular. The $B$'s are also  regular and therefore such is $N$.
This concludes the proof of Theorem~\ref{t1.1}.

\hskip13cm $\Box$

\end{document}